\documentclass[12pt]{amsart}
\usepackage{times}
\usepackage{amsfonts}
\usepackage{amsthm}
\usepackage{amsmath}
\advance \topmargin by -\headheight
\advance \topmargin by -\headsep
\evensidemargin \oddsidemargin
\newtheorem{theorem}{Theorem}[section]

\newtheorem{lemma}[theorem]{Lemma}
\newtheorem{corollary}[theorem]{Corollary}

\newtheorem{remark}[theorem]{Remark}

\begin{document}

\title{ A 
propagation 
property of free entropy dimension }

\begin{abstract} Let $M$ be a tracial von Neumann algebra and $A$ be a weakly dense unital  $C^*$-subalgebra of $M$.  We say that a set $X$ is a $W^*$-generating set for $M$ if the von Neumann algebra generated by $X$ is $M$ and that $X$ is a $C^*$-generating set for $A$ if the unital $C^*$-algebra generated by $X$ is $A$.  For any finite $W^*$-generating set $X$ for $M$ we show that $\delta_0(X) \leq \sup \{\delta_0(Y): \text{$Y$ is a finite $C^*$-generating set for $A$} \}$. It follows that if $\sup 
\{\delta_0(Y): \text{$Y$ is a finite $C^*$-generating set for $C^*_{red}(\mathbb F_2)$}\} < \infty$, then the free group 
factors are all nonisomorphic.  \end{abstract}

\author{Kenley Jung}

\address{Department of Mathematics, University of California,
Los Angeles, CA 90095-1555,USA}

\email{kjung@math.ucla.edu}
\thanks{Research supported in part by the NSF}
\maketitle

\section*{Introduction} 

The past decade has seen the introduction of several notions of dimension for tracial von Neumann algebras defined in terms of generating sets.  Given a finite set $X$ in a tracial von Neumann algebra, Voiculescu defined in \cite{dvv:entropy2} and \cite{dvv:entropy3} the free entropy dimension of $X$, $\delta(X)$, and the modified free entropy dimension of $X$, $\delta_0(X)$.  These notions of dimension were defined in terms of a "microstate free entropy".  There is another such quantity using Hausdorff outer measure on the microstate spaces, denoted by $\mathbb H(X)$ introduced in \cite{jung:fedm}.  Voiculescu went further and developed (via a route which bypasses microstates and uses free difference quotients) another type of entropy, sometimes called the "non-microstates entropy".  This entropy provides a way to define a "non-microstates free entropy dimension" for $X$, denoted by $\delta^*(X)$.  Computations with $\delta^*(X)$ gave rise to yet another notion of dimension, discussed in \cite{Aagard} and  \cite{connes-shlyakht:l2betti}, denoted by $\delta^{\star}(X)$.  Connes and Shlyakhtenko also introduced in \cite{connes-shlyakht:l2betti} a notion of dimension for $X$ denoted by $\Delta(X)$, defined in terms of module dimension.  The relationship between all these quantities is as follows:

\[ \mathbb H(X) \leq \delta_0(X) \leq \delta(X) \leq \delta^{\star}(X) \leq \delta^{\ast}(X) \leq \Delta(X).\]    

\noindent It is not known whether strict inequality can hold between any of the two quantities above. 

A fundamental question (first raised by Voiculescu in \cite{dvv:entropysurvey} for $\delta_0$) which has generated much of the interest in these notions of dimension is whether any of them are von Neumann algebra invariants, i.e., if $X$ and $Y$ are finite sets in $(M,\varphi)$ and they generate the same von Neumann algebra, then is it always the case that $\dim (X) = \dim (Y)$ where $\dim \in \{\mathbb H, \delta_0, \delta, \delta^{\star}, \delta^*, \Delta\}$? 

The first three microstates dimensions are "global" in flavor: one must make asymptotic estimates on the fractal dimension of Euclidean subsets that model the structure of the von Neumann algebra.  Although this can be difficult, it lead to the solutions of several longstanding problems (\cite{dvv:entropy3}, \cite{ge}). The three dominating terms have a "local" quality where derivations, bimodules, and Murray-von Neumann dimension come into play; while certain properties enjoyed by the microstate quantities are lost, one gains considerable flexibility: it was shown in \cite{connes-shlyakht:l2betti} and \cite{mineyev-shlyakht} that if $X$ consists of a finite set of group elements, then $\delta^{\star}(X) = \delta^*(X) = \Delta(X) = \beta_2^1(G) - \beta_2^0(G) +1$ where $G$ is the group generated by $X$ and $\beta_2^i(G)$ are the $L^2$-Betti numbers of the group $G$.  

While it seems overly optimistic to hope for equality among any of the 
terms in full generality, one can play the safer game of taking known 
properties for one dimension and trying to prove the same property for 
another dimension (with some extra hypothesis placed on the von Neumann 
algebra, if necessary).  The remark we want to make here is an example of 
this.

In \cite{thom}, Thom examines some of the notions introduced in 
\cite{connes-shlyakht:l2betti} and reduces some von Neumann algebra 
difficulties to $C^*$-algebra ones.  One of the main tools used in this 
reduction is a noncommutative generalization of Lusin's Theorem 
(\cite{tomita}, \cite{saito}, \cite{pedersen}), henceforth simply referred 
to as Lusin's Theorem.  The 
observation we make here is that Lusin's Theorem has a similar 
consequence in the context of microstate free entropy dimension.  We will 
show that if $A$ is a weakly dense $C^*$-algebra of $M$ and 
$\overline{\delta_0}(A)$ is the supremum over all values of the form 
$\delta_0(F)$ where $F$ is a finite $C^*$-generating set for $A$ (the unital $C^*$-algebra generated by $F$ is $A$), then for any finite 
$W^*$-generating set $X$ of $M$ (the von Neumann algebra generated by $X$ is 
$M$), $\delta_0(X) \leq \overline{\delta_0}(A)$.  In 
particular, $C^*$-algebra invariance for $\delta_0$ propagates to von Neumann algebra invariance for $\delta_0$ in the following sense: if for 
any weakly dense $C^*$-subalgebra $A$ of $M$ and any two finite $C^*$-generating 
sets $X$ and $Y$ for $A$, $\delta_0(X)=\delta_0(Y)$, then $\delta_0$ is a 
von Neumann algebra invariant for $M$.  By \cite{ken} and \cite{florin}, 
it also follows that if $\overline{\delta_0}(C^*_{red}(\mathbb F_2)) < 
\infty$, then the free group factors are nonisomorphic.  The proof of the 
main inequality is short and trivial once one has Lusin's Theorem.

\section{The propagation property for $\delta_0$}

Assume throughout that $(M,\varphi)$ is a tracial von 
Neumann algebra and $A$ is a weakly dense C*-subalgebra of $M$.  For a 
given set of elements $X$, we denote by $\Gamma(X;m,k,\gamma)$ the 
nonselfadjoint microstates of $X$ with parameters $(m,k,\gamma)$, i.e.,
$\Gamma(X;m,k,\gamma)$ consists of all $n$-tuples $(a_1,\ldots, a_n)$ of 
$k \times k$ complex matrices 
such that for any $1 \leq p \leq m$, $1\leq i_1,\ldots, i_p \leq n$, 
and 
$j_1,\ldots, j_p \in \{1,*\}$,

\[ |\varphi(x_{i_1}^{j_1} \cdots x_{i_p}^{j_p}) - tr_k(a_{i_1}^{j_1} 
\cdots a_{i_p}^{j_p})| < \gamma.\]

\noindent We can regard $\Gamma(X:m,k,\gamma)$ as metric spaces with 
respect to the inner product metric generated by $|(a_1,\ldots, a_n)|_2 = 
\sum_{j=1}^n tr_k(a_j^*a_j)$.  For any $\epsilon >0$ and metric space 
$\Omega$, denote by $K_{\epsilon}(\Omega)$ the minimum number of open 
$\epsilon$-balls required to cover $\Omega$.  One can define successively,

\[ \mathbb K_{\epsilon}(X;m,\gamma) = \limsup_{k \rightarrow \infty} \log 
k^{-2} \cdot \log(K_{\epsilon}(\Gamma(X;m,k,\gamma))),
\] 

\[ \mathbb K_{\epsilon}(X) = \inf \{K_{\epsilon}(X;m,\gamma): m \in 
\mathbb N, \gamma >0\}.
\]

\noindent One finally defines the (modified) free entropy dimension of $X$ 
by

\[ \delta_0(X) = \limsup_{\epsilon \rightarrow 0} \frac{\mathbb 
K_{\epsilon}(X)}{|\log 
\epsilon|}.
\]

\noindent See \cite{djs} for a further discussion of this nonselfadjoint 
adaptation of these quantities.   We will sometimes write $\mathbb K_{\epsilon}(x_1, \ldots, x_n: \cdot)$ and $\delta_0(x_1,\ldots, x_n)$ for the quantities $K_{\epsilon}(X: \cdot)$, and $\delta_0(X)$ where $X = \{x_1,\ldots, x_n\}$.  

If $A$ is finitely generated then define $\overline{\delta_0}(A)$ to be the supremum over $\delta_0(F)$ 
where $F$ is a finite $C^*$-generating set for $A$; otherwise, set $\overline{\delta_0}(A) = \infty$.

We will need one prepatory result before the theorem.  

\begin{lemma} If $x, p \in M$ with $p$ a projection, then 
$\delta_0(xp) \leq 4 \varphi(p) - 2\varphi(p)^2$.  \end{lemma}

\begin{proof} Fix $R > \|x\| +1$ and denote by $e$ the projection onto the 
range of $xp$. Observe now that $\delta_0(xp) \leq \delta_0(xp, e) = 
\delta_0(e, xe, e)$.  Pick a sequence of projections $\langle e_k 
\rangle_{k=1}^{\infty}$ such that the rank of $e_k$ is $[\varphi(e)k]$.  
If we consider the relative microstate quantities $\Xi(xe, e;m,k,\gamma)$ 
with respect 
to this sequence then we have by \cite{jung:hid} that

\begin{eqnarray} \delta_0(xe, e) = \delta_0(e) + \limsup_{\epsilon 
\rightarrow 0} \frac {\mathbb K_{\epsilon}(\Xi(xe, e))}{|\log \epsilon|}. 
\end{eqnarray}

\noindent Technically speaking, this equation works only for selfadjoint 
microstate spaces, but the proof can be easily adapted to the 
nonselfadjoint situation verbatim to arrive at the equation above (passage 
from the selfadjoint to the nonselfadjoint setting in cases involving 
$\delta_0$ is for the most part trivial).  To prove the inequality it thus 
suffices to bound the second summand above.  Towards this end, pick an 
$\epsilon$-net $\langle \xi_{jk} \rangle_{j \in \Omega_k}$ for the 
$R$-ball of $M_k(\mathbb C)e_k$ (with respect to the $| \cdot |_2$-norm) 
with $\# \Omega_k < \left (\frac{2R}{\epsilon}\right)^{2k^2 tr_k(e_k)}$ 
(by an $\epsilon$-cover we mean a set whose $\epsilon$-neighborhood 
contains the set in questions).  Choose $m \in \mathbb N$ and $\gamma >0$ 
appropriately so that if $(a, b) \in \Xi(xe,e;m,k,\gamma)$ then $|ae_k - 
a|_2 < \epsilon$, $|ae_k|_2 <R$, $|b - e_k|_2 < \epsilon$. It then follows 
that $\langle (\xi_{jk}e_k, e_k) \rangle_{j \in S_k}$ is a $3\epsilon$ 
-net for $\Xi(xe,e)$ with respect to the $| \cdot |_2$ norm.  
To see this just observe that there
exists a $j_0 \in \Omega_k$ such that $| ae_k - \xi_{j_0k}|_2 < 
\epsilon$, whence $|a - \xi_{j_0k}|_2 < 2 \epsilon$.  Thus, $|(a, b)
- (\xi_{j_0k}, e_k)|_2 < 3 \epsilon$.

We now have the bound

\begin{eqnarray*} \limsup_{\epsilon \rightarrow 0} \mathbb K_{4 
\epsilon}(\Xi(xe,e)) & \leq & \limsup_{k \rightarrow \infty} k^{-2} \cdot 
\log(K_{3\epsilon}(\Xi(xe,e;m,k,\gamma))) \\ & \leq & \limsup_{k 
\rightarrow \infty} k^{-2} \cdot 2k^2 tr_k(e_k) \log 
\left(\frac{2R}{\epsilon} 
\right) \\ & \leq & 2 \varphi(e) (|\log \epsilon| + \log (2R)).\\ 
\end{eqnarray*}

\noindent From this and \cite{dvv:entropy3} we have that $(1)$ is dominated by

\begin{eqnarray*} \delta_0(e) + 2 \varphi(e) & = & 1 - \varphi(e)^2 - (1 - 
\varphi(e))^2 + 2 \varphi(e) \\ & = & 4 \varphi(e) - 2 \varphi(e)^2 \\ & 
\leq & 4 \varphi(p) - 2\varphi(p)^2.\end{eqnarray*} \end{proof}

\begin{theorem} If $X$ is a finite $W^*$-generating set for $M$, 
then $\delta_0(X) \leq \overline{\delta_0}(A)$. \end{theorem}

\begin{proof} Suppose $X = \{x_1,\ldots, x_n\}$ and $\epsilon >0$.  By 
Lusin's Theorem (see \cite{tomita}, \cite{saito}, \cite{pedersen}) we can 
find for each $1\leq i \leq n$ a $y_i \in A$ and a 
projection $p_i \in M$ such that $x_ip_i = y_ip_i$ and $\varphi(p_i) > 1- 
\epsilon$.  Set $Y = \{y_1,\ldots, y_n\}$. Using the fact that $\delta_0$ 
is an algebraic invariant and monotonic with respect to adding finite sets 
already in the von Neumann algebra of the given set 
(\cite{dvv:strengthened}),

\begin{eqnarray*} \delta_0(X) & \leq & \delta_0(X \cup \{p_1,\ldots, p_n\} \cup G \cup Y)\\ & \leq & \delta_0( \{x_1p_1,\ldots, x_np_n\} \cup \{x_1p_1^{\bot},\ldots, x_n p_n^{\bot}\} \cup \{p_1,\ldots, p_n\} \cup G \cup Y) \\ & = & \delta_0(\{x_1p_1^{\bot},\ldots, x_n p_n^{\bot}\} \cup \{p_1,\ldots, p_n\} \cup G \cup Y) \\ & \leq & \Sigma_{j=1}^n \delta_0(x_j p_j^{\bot}) + \Sigma_{j=1}^n \delta_0(p_j) + \delta_0(G \cup Y) \\ & \leq & n(4\epsilon -2 \epsilon^2) + n(2\epsilon -\epsilon^2) + \overline{\delta_0}(A). 
\end{eqnarray*}

\noindent Since $\epsilon$ was independent of $n$ we have the desired 
inequality. \end{proof}

\begin{corollary} If $A_1$ and $A_2$ are weakly dense C*-subalgebras of 
$M$, then $\overline{\delta_0}(A_1) = \overline{\delta_0}(A_2)$. 
\end{corollary}

For $\alpha >0$, denote by $M_{\alpha}$ the amplification of $M$ by 
$\alpha$.  $\mathcal F(M)$ will be the fundamental group of $M$, i.e., 
$\mathcal F(M) = \{ \alpha \in (0,\infty]: M_{\alpha} \simeq M\}$.  We 
have the following easy consequence of Theorem 1.2:

\begin{corollary} If $1 < \overline{\delta_0}(A) < \infty$, then $\mathcal F(M) = \{1\}$.
\end{corollary}

\begin{proof} As $\mathcal F(M)$ is a multiplicative subgroup of the nonzero positive reals, it suffices to show that if $1 < \beta$, then $\beta \notin \mathcal F(M)$.  So, suppose $1 < \beta$.  $M$ is finitely generated so $M_{\beta}$ is as well.  Fix a finite generating set $X$ for $M_{\beta}$.  Then by the assumption, the lemma above, and Lemma 5.2 of \cite{jung:index},

\[ \delta_0(X) \leq 1 - \beta^{-2} + \beta^{-2} \cdot  \overline{\delta_0(A)} < \overline{\delta_0(A)}. \]

\noindent This being true for any such set $X$, it follows that $M_{\beta}$ cannot be isomorphic to $M$, i.e., $\beta \notin \mathcal F(M)$.  This completes the proof.
\end{proof}

\begin{remark} Notice that if $M$ is diffuse and embeddable into $\mathcal R^{\omega}$, then it is automatic from \cite{jung:hf} then $1 \leq \overline{\delta_0}(A)$. 
\end{remark}

Using the compression formulae from \cite{florin} followed by the fact that the free group factors are either all isomorphic or all non-isomorphic (\cite{ken}) we get: 

\begin{corollary} If for some $n \geq 1$, $M = L(\mathbb F_n)$ and 
$\overline{\delta_0(A)} < \infty$, then the free group factors are 
non-isomorphic. \end{corollary}

\noindent{\it Acknowledgments.} I'd like to thank Dan Voiculescu and Dimitri Shlyakhtenko for suggesting that I write up this observation.

\bibliographystyle{amsplain}

\providecommand{\bysame}{\leavevmode\hbox to3em{\hrulefill}\thinspace}

\end{document}